\documentclass[12pt]{elsart}
\usepackage{amsmath,amsfonts,amssymb,amsbsy}
\usepackage[cp1251]{inputenc}
\usepackage[english]{babel}
\usepackage{color}
\usepackage[pdftex]{graphicx}

\begin{document}

\begin{frontmatter}
\title{
An asymptotic model for a thin bonded elastic layer coated with an elastic membrane}
\author{I.~Argatov},
\corauth[cor]{Corresponding author.}
\author{G.~Mishuris\corauthref{cor}}
\ead{ggm@aber.ac.uk}
\address{Institute of Mathematics and Physics, Aberystwyth University,
Ceredigion SY23 3BZ, Wales, UK}
\begin{abstract}
The deformation problem for a transversely isotropic elastic layer bonded to a rigid substrate and coated with a very thin elastic layer made of another transversely isotropic material is considered.
The leading-order asymptotic models (for compressible and incompressible layers) are constructed based on the simplifying assumptions that the generalized plane stress conditions apply to the coating layer, and the flexural stiffness of the coating layer is negligible compared to its tensile stiffness.
\end{abstract}

\begin{keyword}
Deformation problem \sep thin layer \sep elastic coating \sep transversely isotropic \sep asymptotic solution
\end{keyword}
\end{frontmatter}

\setcounter{equation}{0}

\section{Introduction}

Some natural biological tissues such as articular cartilage possess an inhomogeneous, layered structure with anisotropic material properties.
In particular, morphological studies of adult articular cartilage \cite{Clarke1971,Hughes_et_al2005} show three different zones of preferred collagen fiber bundle orientation. The superficial zone formed by tangentially oriented collagen fibrils provides a thin layer with a high tensile stiffness in the direction parallel to articular surface.
It was shown \cite{Korhonen2002} that the high transverse stiffness of the superficial tissue layer (characterized by tangentially oriented collagen fibrils) is important in controlling the deformation response of articular cartilage.
Generally speaking, the surface layer in biomaterials usually has different mechanical properties than the underlying bulk material.
Due to this circumstance, the mechanical deformation behavior is strongly influenced by the complex interaction of these layers \cite{RahmanNewaz2000}.

In the present paper, we consider the deformation problem for a transversely isotropic elastic layer reinforced with a thin elastic membrane ideally attached to one surface, while the other surface is bonded to a rigid substrate.
Following \cite{AlexandrovMkhitaryan1985,RahmanNewaz1997}, it is assumed that the reinforcing layer is very thin (with respect to a characteristic size of the applied load) so that its deformation can be treated in the framework of the generalized plane stress state. Moreover, it is assumed that the flexural stiffness of the coating layer is negligible compared to its tensile stiffness.
Thus, the reinforcing thin layer is regarded as an elastic membrane.

The rest of the paper is organized as follows. In Section~\ref{sec:44.1.1}, we formulate the three-dimensional boundary conditions for a coated elastic layer.
The deformation problem formulation itself is given in Section~\ref{sec:44.1.2}.
Asymptotic analysis of the deformation problem is presented in Section~\ref{sec:44.1.3}.
The main result of the present paper is presented by the leading-order asymptotic models for the local indentation of the coated elastic layer developed in Section~\ref{sec:44.1.4} for the cases of compressible and incompressible layer.
Finally, Section~\ref{sec:44.1.4D} contains some discussion of the obtained asymptotic models and outlines our conclusion.

\section{Boundary conditions for a coated elastic layer}
\label{sec:44.1.1}

\begin{figure}[h!]
    \centering
\vbox{
    \includegraphics [scale=0.5]{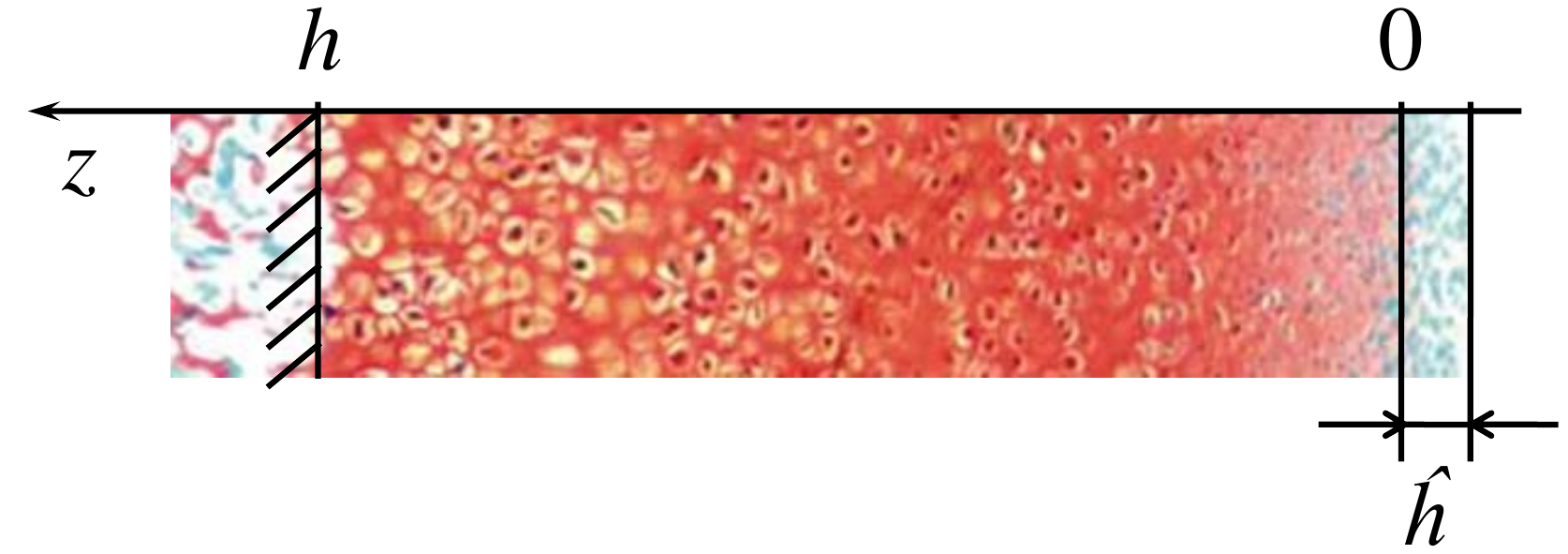}
    \caption{An elastic coated layer as a model for articular cartilage (the histological image of articular cartilage is taken from the paper \cite{PatelMao2003}). The regions $-\hat{h} \le z\le 0$ and $z\geq h$ represent the superficial zone and the subchondral bone, respectively.}
    \label{figureartC.pdf}
   }
\end{figure}

We consider a very thin transversely isotropic elastic coating layer (of uniform thickness $\hat{h}$) ideally attached to an elastic layer (of thickness $h$) made of another transversely isotropic material (see Fig.~\ref{figureartC.pdf}).
Let the five independent elastic constants of the elastic layer and its coating are denoted by
$A_{11}$, $A_{12}$, $A_{13}$, $A_{33}$, $A_{44}$
and $\hat{A}_{11}$, $\hat{A}_{12}$, $\hat{A}_{13}$, $\hat{A}_{33}$, $\hat{A}_{44}$, respectively.

Under the assumption that the two layers are in perfect contact with one another along their common interface, $z=0$, the following boundary conditions of continuity (interface conditions of perfect bonding) should be satisfied:
\begin{equation}
\hat{\bf v}({\bf y},0)={\bf v}({\bf y},0),\quad
\hat{w}({\bf y},0)=w({\bf y},0),
\label{44ch(4.1)}
\end{equation}
\begin{equation}
\hat{\sigma}_{3j}({\bf y},0)=\sigma_{3j}({\bf y},0),\quad j=1,2,3.
\label{44ch(4.2)}
\end{equation}
Here, $(\hat{\bf v},\hat{w})$ is the displacement vector of the elastic coating layer,
$\hat{\sigma}_{ij}$ are the corresponding components of stress.
In what follows, we make use of the Cartesian coordinate system $({\bf y},z)$, where ${\bf y}=(y_1,y_2)$ are the in-plane coordinates.

On the upper surface of the two-layer system, $z=-\hat{h}$, we impose the boundary conditions of normal loading with no tangential tractions
\begin{equation}
\hat{\sigma}_{31}({\bf y},-\hat{h})=\sigma_{32}({\bf y},-\hat{h})=0,\quad
\hat{\sigma}_{33}({\bf y},-\hat{h})=-p({\bf y}),
\label{44ch(4.3)}
\end{equation}
where $p({\bf y})$ is a specified function of external loading.

Following Rahman and Newaz \cite{RahmanNewaz2000}, we simplify the deformation analysis of the elastic coating layer based on the following two assumptions: (1) the coating layer is assumed to be very thin, so that the generalized plane stress conditions apply; (2)~the flexural stiffness of the coating layer in the $z$-direction is negligible compared to its tensile stiffness.

In the absence of body forces, the equilibrium equations for an infinitesimal element of the coating layer are
\begin{equation}
\frac{\partial\hat{\sigma}_{j1}}{\partial y_1}
+\frac{\partial\hat{\sigma}_{j2}}{\partial y_2}
+\frac{\partial\hat{\sigma}_{j3}}{\partial z}=0,\quad
j=1,2,3.
\label{44ch(4.4)}
\end{equation}

The stress-strain relationship for the transversely isotropic elastic coating layer is given by
\begin{equation}
\left(
\begin{array}{c}
\hat{\sigma}_{11} \\
\hat{\sigma}_{22} \\
\hat{\sigma}_{33} \\
\hat{\sigma}_{23} \\
\hat{\sigma}_{13} \\
\hat{\sigma}_{12}
\end{array}
\right)=
\left[
\begin{array}{cccccc}
\hat{A}_{11} & \hat{A}_{12} & \hat{A}_{13} & 0 & 0 & 0 \\
\hat{A}_{12} & \hat{A}_{11} & \hat{A}_{13} & 0 & 0 & 0 \\
\hat{A}_{13} & \hat{A}_{13} & \hat{A}_{33} & 0 & 0 & 0 \\
0 & 0 & 0 & 2\hat{A}_{44} & 0 & 0 \\
0 & 0 & 0 & 0 & 2\hat{A}_{44} & 0 \\
0 & 0 & 0 & 0 & 0 & 2\hat{A}_{66}
\end{array}
\right]
\left(
\begin{array}{c}
\hat{\varepsilon}_{11} \\
\hat{\varepsilon}_{22} \\
\hat{\varepsilon}_{33} \\
\hat{\varepsilon}_{23} \\
\hat{\varepsilon}_{13} \\
\hat{\varepsilon}_{12}
\end{array}
\right),
\label{44ch(4.6)}
\end{equation}
where $2\hat{A}_{66}=\hat{A}_{11}-\hat{A}_{12}$.

Integrating Eqs.~(\ref{44ch(4.4)}) through the thickness of coating layer and taking into account the interface and boundary conditions (\ref{44ch(4.2)}) and (\ref{44ch(4.3)}), we get
\begin{equation}
\hat{h}\biggl(
\frac{\partial\hat{\bar{\sigma}}_{j1}}{\partial y_1}
+\frac{\partial\hat{\bar{\sigma}}_{j2}}{\partial y_2}
\biggr)=-\sigma_{j3}\bigr\vert_{z=0},\quad j=1,2,
\label{44ch(4.7)}
\end{equation}
\begin{equation}
\hat{h}\biggl(
\frac{\partial\hat{\bar{\sigma}}_{13}}{\partial y_1}
+\frac{\partial\hat{\bar{\sigma}}_{23}}{\partial y_2}
\biggr)=-\sigma_{33}\bigr\vert_{z=0}-p.
\label{44ch(4.8)}
\end{equation}
Here, $\hat{\bar{\sigma}}_{ij}$ are the averaged stresses, i.e.,
$$
\hat{\bar{\sigma}}_{ij}({\bf y})=\frac{1}{\hat{h}}
\int\limits_{-\hat{h}}^0\hat{\sigma}_{ij}({\bf y},z)\,dz.
$$

Under the simplifying assumptions made above, we have
\begin{equation}
\hat{\bar{\sigma}}_{13}=\hat{\bar{\sigma}}_{23}=\hat{\bar{\sigma}}_{33}=0.
\label{44ch(4.9)}
\end{equation}

Hence, Eq.~(\ref{44ch(4.8)}) immediately implies that
\begin{equation}
\sigma_{33}\bigr\vert_{z=0}=-p.
\label{44ch(4.10)}
\end{equation}

Moreover, in view of (\ref{44ch(4.9)}), the averaged strain
$\hat{\bar{\varepsilon}}_{33}$ must satisfy the equation
$$
\hat{A}_{13}\hat{\bar{\varepsilon}}_{11}+\hat{A}_{13}\hat{\bar{\varepsilon}}_{22}
+\hat{A}_{33}\hat{\bar{\varepsilon}}_{33}=0.
$$

Therefore, the in-plane averaged stress-strain relationship takes the form
\begin{equation}
\left(
\begin{array}{c}
\hat{\bar{\sigma}}_{11} \\
\hat{\bar{\sigma}}_{22} \\
\hat{\bar{\sigma}}_{12}
\end{array}
\right)=
\left[
\begin{array}{ccc}
\hat{\bar{A}}_{11} & \hat{\bar{A}}_{12} & 0 \\
\hat{\bar{A}}_{12} & \hat{\bar{A}}_{11} & 0 \\
0 & 0 & 2\hat{\bar{A}}_{66}
\end{array}
\right]
\left(
\begin{array}{c}
\hat{\bar{\varepsilon}}_{11} \\
\hat{\bar{\varepsilon}}_{22} \\
\hat{\bar{\varepsilon}}_{12}
\end{array}
\right),
\label{44ch(4.11)}
\end{equation}
where we introduced the notation
\begin{equation}
\hat{\bar{A}}_{11}=\hat{A}_{11}-\frac{\hat{A}_{13}^2}{\hat{A}_{33}},\quad
\hat{\bar{A}}_{12}=\hat{A}_{12}-\frac{\hat{A}_{13}^2}{\hat{A}_{33}},\quad
2\hat{\bar{A}}_{66}=\hat{\bar{A}}_{11}-\hat{\bar{A}}_{12}.
\label{44ch(4.12)}
\end{equation}

On the other hand, in view of the interface conditions (\ref{44ch(4.1)}), we will have
\begin{equation}
\hat{\bar{\varepsilon}}_{11}=\varepsilon_{11}\bigr\vert_{z=0},\quad
\hat{\bar{\varepsilon}}_{22}=\varepsilon_{22}\bigr\vert_{z=0},\quad
\hat{\bar{\varepsilon}}_{12}=\varepsilon_{12}\bigr\vert_{z=0},
\label{44ch(4.13)}
\end{equation}
where $\varepsilon_{11}$, $\varepsilon_{22}$, and $\varepsilon_{12}$ are the in-plane strains in the coated elastic layer $z\in(0,h)$.

Therefore, taking Eqs.~(\ref{44ch(4.11)}) and (\ref{44ch(4.13)}) into account, we transform the boundary conditions (\ref{44ch(4.7)}) as follows:
$$
-\frac{1}{\hat{h}}\sigma_{31}\bigr\vert_{z=0}=
\frac{\partial}{\partial y_1}\biggl(
\hat{\bar{A}}_{11}\frac{\partial v_1}{\partial y_1}
+\hat{\bar{A}}_{12}\frac{\partial v_2}{\partial y_2}\biggr)
+\hat{\bar{A}}_{66}\frac{\partial}{\partial y_2}\biggl(
\frac{\partial v_1}{\partial y_2}+\frac{\partial v_2}{\partial y_1}
\biggr),
$$
$$
-\frac{1}{\hat{h}}\sigma_{32}\bigr\vert_{z=0}=
\hat{\bar{A}}_{66}\frac{\partial}{\partial y_1}\biggl(
\frac{\partial v_1}{\partial y_2}+\frac{\partial v_2}{\partial y_1}
\biggr)
+\frac{\partial}{\partial y_2}\biggl(
\hat{\bar{A}}_{12}\frac{\partial v_1}{\partial y_1}
+\hat{\bar{A}}_{11}\frac{\partial v_2}{\partial y_2}\biggr).
$$

Finally, the above boundary conditions can be rewritten in the matrix form as
\begin{equation}
\sigma_{31}{\bf e}_1+\sigma_{32}{\bf e}_2\bigr\vert_{z=0}
=-\hat{\mathfrak{L}}(\nabla_y){\bf v}\bigr\vert_{z=0},
\label{44ch(4.14)}
\end{equation}
where $\hat{\mathfrak{L}}(\nabla_y)$ is a $2\times 2$ matrix differential operator such that
\begin{equation}
\begin{array}{rcl}
\hat{\mathfrak{L}}_{\alpha\alpha}(\nabla_y) & = & \displaystyle
\hat{h}\hat{\bar{A}}_{11}\frac{\partial^2}{\partial y_\alpha^2}
+\hat{h}\hat{\bar{A}}_{66}\frac{\partial^2}{\partial y_{3-\alpha}^2}
\phantom{{}_{\Bigr)}}\\
\hat{\mathfrak{L}}_{\alpha\beta}(\nabla_y) & = & \displaystyle
\hat{h}\bigl(\hat{\bar{A}}_{12}+\hat{\bar{A}}_{66}\bigr)
\frac{\partial^2}{\partial y_\alpha\partial y_\beta},\quad
\alpha,\beta=1,2,\quad \alpha\not=\beta.
\end{array}
\label{44ch(4.14L)}
\end{equation}

Thus, the deformation problem for an elastic layer coated with a very thin flexible elastic layer is reduced to that for the elastic layer without coating, but subjected to a different set of boundary conditions (\ref{44ch(4.10)}) and (\ref{44ch(4.14)}) on the surface $z=0$.

\section{Deformation problem formulation}
\label{sec:44.1.2}

Now, let us consider a relatively thin transversely isotropic elastic layer of uniform thickness, $h$, coated with an infinitesimally thin elastic membrane and bonded to a rigid substrate, so that
\begin{equation}
{\bf v}\bigr\vert_{z=h}={\bf 0},\quad w\bigr\vert_{z=h}=0.
\label{44ch(4.16)}
\end{equation}

In the absence of body forces, the vector $({\bf v},w)$ of displacements in the elastic layer satisfies the Lam\'e system
\begin{equation}
\begin{array}{rcl}
\displaystyle
A_{66}\Delta_y{\bf v}+(A_{11}-A_{66})\nabla_y\nabla_y\!\cdot{\bf v}+A_{44}
\frac{\partial^2{\bf v}}{\partial z^2}+(A_{13}+A_{44})\frac{\partial}{\partial z}\nabla_y w & = & {\bf 0},\phantom{{}_{\Bigr)}}\\
\displaystyle
A_{44}\Delta_y w+A_{33}\frac{\partial^2 w}{\partial z^2}+(A_{13}+A_{44})\frac{\partial}{\partial z}\nabla_y\!\cdot{\bf v} & = & 0.
\end{array}
\label{44ch(4.17)}
\end{equation}

Assuming that the coated layer is loaded by a normal load and denoting the load density by $p$, we require that
\begin{equation}
\sigma_{33}\bigr\vert_{z=0}=-p.
\label{44ch(4.18)}
\end{equation}

Based on the analysis performed in Section~\ref{sec:44.1.1}, the influence of the elastic membrane (coating layer) is introduced by the boundary condition
\begin{equation}
\sigma_{31}{\bf e}_1+\sigma_{32}{\bf e}_2\bigr\vert_{z=0}
=-\hat{\mathfrak{L}}(\nabla_y){\bf v}\bigr\vert_{z=0},
\label{44ch(4.19)}
\end{equation}
where $\hat{\mathfrak{L}}(\nabla_y)$ is the matrix differential operator defined by formulas (\ref{44ch(4.14L)}).

Taking into account the stress-strain relationship
$$
\left(
\begin{array}{c}
\sigma_{11} \\
\sigma_{22} \\
\sigma_{33} \\
\sigma_{23} \\
\sigma_{13} \\
\sigma_{12}
\end{array}
\right)=
\left[
\begin{array}{cccccc}
A_{11} & A_{12} & A_{13} & 0 & 0 & 0 \\
A_{12} & A_{11} & A_{13} & 0 & 0 & 0 \\
A_{13} & A_{13} & A_{33} & 0 & 0 & 0 \\
0 & 0 & 0 & 2A_{44} & 0 & 0 \\
0 & 0 & 0 & 0 & 2A_{44} & 0 \\
0 & 0 & 0 & 0 & 0 & 2A_{66}
\end{array}
\right]
\left(
\begin{array}{c}
\varepsilon_{11} \\
\varepsilon_{22} \\
\varepsilon_{33} \\
\varepsilon_{23} \\
\varepsilon_{13} \\
\varepsilon_{12}
\end{array}
\right),
$$
we rewrite Eqs.~(\ref{44ch(4.18)}), (\ref{44ch(4.19)}) as follows:
\begin{equation}
A_{13}\frac{\partial v_1}{\partial y_1}+A_{13}\frac{\partial v_2}{\partial y_2}
+A_{33}\frac{\partial w}{\partial z}\biggr\vert_{z=0}=-p,
\label{44ch(4.20)}
\end{equation}
\begin{equation}
A_{44}\Bigl(\nabla_y w+\frac{\partial{\bf v}}{\partial z}\Bigr)\Bigr\vert_{z=0}
=-\hat{\mathfrak{L}}(\nabla_y){\bf v}\bigr\vert_{z=0}.
\label{44ch(4.21)}
\end{equation}

Equations (\ref{44ch(4.16)}), (\ref{44ch(4.17)}), (\ref{44ch(4.20)}), and (\ref{44ch(4.21)}) comprise the deformation problem for the coated transversely isotropic elastic layer bonded to a rigid substrate.

\section{Asymptotic analysis of the deformation problem}
\label{sec:44.1.3}

Let $h_*$ be a characteristic length of the external load distribution. Denoting by $\varepsilon$ a small positive parameter, we require that
\begin{equation}
h=\varepsilon h_*
\label{44ch(4.22)}
\end{equation}
and introduce the so-called ``stretched'' dimensionless normal coordinate
$$
\zeta=\frac{z}{\varepsilon h_*}.
$$
In addition, we non-dimensionalize the in-plane coordinates by the formulas
$$
\eta_i=\frac{y_i}{h_*}, \quad i=1,2,\quad \mbox{\boldmath $\eta$}=(\eta_1,\eta_2),
$$
so that
$$
\frac{\partial }{\partial z}=\frac{1}{\varepsilon h_*}\frac{\partial }{\partial\zeta},
\quad \nabla_y=\frac{1}{h_*}\nabla_\eta.
$$

Moreover, we assume that the tensile stiffness of the coating layer is relatively high, i.e., $\hat{\bar{A}}_{11}>\!\!> A_{11}$, and so on. To fix our ideas, we consider the situation when
\begin{equation}
\hat{\mathfrak{L}}(\nabla_y)=\varepsilon^{-1}
\hat{\mathfrak{L}}^*(\nabla_y),
\label{44ch(4.2L)}
\end{equation}
that is, in particular, the ratio $A_{11}/\hat{\bar{A}}_{11}$ is of the order of $\varepsilon$.

Employing the perturbation algorithm \cite{Gol'denveizer1962-77}, we represent the solution to the deformation problem (\ref{44ch(4.16)}), (\ref{44ch(4.17)}), (\ref{44ch(4.20)}), (\ref{44ch(4.21)}) as follows:
\begin{eqnarray}
{\bf v} & = & \varepsilon^2{\bf v}^1(\mbox{\boldmath $\eta$},\zeta)
+\ldots,
\label{44ch(4.23)}\\
 w & = & \varepsilon w^0(\mbox{\boldmath $\eta$},\zeta)
+\varepsilon^3 w^2(\mbox{\boldmath $\eta$},\zeta)+\ldots\,.
\label{44ch(4.24)}
\end{eqnarray}
Here, for the sake of brevity, we include only non-vanishing terms.

It can be shown (see, in particular, \cite{Argatov2004thin}) that the leading-order term in (\ref{44ch(4.24)}) is given by
\begin{equation}
w^0(\mbox{\boldmath $\eta$},\zeta)=\frac{h_* p}{A_{33}}(1-\zeta),
\label{44ch(4.25)}
\end{equation}
whereas the first non-trivial term of the expansion (\ref{44ch(4.23)}) satisfies the problem
$$
A_{44}\frac{\partial^2{\bf v}^1}{\partial \zeta^2}=-(A_{13}+A_{44})\nabla_\eta\frac{\partial w^0}{\partial\zeta},\quad \zeta\in(0,1),
$$
$$
A_{44}\frac{\partial {\bf v}^1}{\partial \zeta}
+\frac{1}{h_*}\hat{\mathfrak{L}}^*(\nabla_\eta){\bf v}^1\biggr\vert_{\zeta=0}
=-A_{44}\nabla_\eta w^0\bigr\vert_{\zeta=0},\quad
{\bf v}^1\bigr\vert_{\zeta=1}={\bf 0}.
$$

Substituting the expansion (\ref{44ch(4.25)}) for $w^0$ into the above equations, we obtain
\begin{equation}
\begin{array}{c}
\displaystyle
A_{44}\frac{\partial^2{\bf v}^1}{\partial \zeta^2}=\frac{A_{13}+A_{44}}{A_{33}}h_*\nabla_\eta p,\quad \zeta\in(0,1),
\\
\displaystyle
A_{44}\frac{\partial {\bf v}^1}{\partial \zeta}
+\frac{1}{h_*}\hat{\mathfrak{L}}^*(\nabla_\eta){\bf v}^1\biggr\vert_{\zeta=0}^{\phantom{\bigr)}}
=-\frac{A_{44}}{A_{33}}h_*\nabla_\eta p,\quad
{\bf v}^1\bigr\vert_{\zeta=1}={\bf 0}.
\end{array}
\label{44ch(4.26)}
\end{equation}

The solution to the boundary-value problem (\ref{44ch(4.26)}) is represented in the form
\begin{equation}
{\bf v}^1=-\frac{A_{13}+A_{44}}{2A_{33}A_{44}}\zeta(1-\zeta)
h_*\nabla_\eta p+(1-\zeta){\bf V}^1(\mbox{\boldmath $\eta$}),
\label{44ch(4.28)}
\end{equation}
where ${\bf V}^1(\mbox{\boldmath $\eta$})$ satisfies the equation
\begin{equation}
\frac{1}{h_*}\hat{\mathfrak{L}}^*(\nabla_\eta){\bf V}^1
-A_{44}{\bf V}^1=\frac{A_{13}-A_{44}}{2A_{33}}h_*\nabla_\eta p
\label{44ch(4.29)}
\end{equation}
on the entire plane $\zeta=0$.

For the second non-trivial term of the expansion (\ref{44ch(4.24)}), we derive the problem
$$
A_{33}\frac{\partial^2 w^2}{\partial \zeta^2}=-(A_{13}+A_{44})\nabla_\eta\!\cdot\frac{\partial {\bf v}^1}{\partial\zeta}-A_{44}\Delta_\eta w^0,\quad \zeta\in(0,1),
$$
$$
A_{33}\frac{\partial w^2}{\partial \zeta}\biggr\vert_{\zeta=0}=-A_{13}\nabla_\eta\!\cdot{\bf v}^1\bigr\vert_{\zeta=0},\quad  w^2\bigr\vert_{\zeta=1}=0.
$$

Now, substituting the expressions (\ref{44ch(4.25)}) and (\ref{44ch(4.28)}) for $w^0$ and ${\bf v}^1$, respectively, into the above equations, we arrive at the problem
\begin{eqnarray}
\frac{\partial^2 w^2}{\partial \zeta^2} & = &
-\bigl[(A_{13}+A_{44})^2(2\zeta-1)+2A_{44}^2(1-\zeta)\bigr]
\frac{h_*\Delta_\eta p}{2A_{33}^2 A_{44}}
\nonumber \\
{ } & { } & {}+\frac{A_{13}+A_{44}}{A_{33}}\nabla_\eta\!\cdot{\bf V}^1,
\quad \zeta\in(0,1),
\label{44ch(4.30)}
\end{eqnarray}
\begin{equation}
\frac{\partial w^2}{\partial \zeta}\biggr\vert_{\zeta=0}=-\frac{A_{13}}{A_{33}}\nabla_\eta\!\cdot{\bf V}^1,\quad  w^2\bigr\vert_{\zeta=1}=0.
\label{44ch(4.31)}
\end{equation}

Integrating Eq.~(\ref{44ch(4.30)}) twice with respect to $\zeta$ and taking into account the boundary condition $(\ref{44ch(4.31)})_2$, we obtain
\begin{eqnarray}
w^2 & = & -\bigl[(A_{13}+A_{44})^2\bigl(2\zeta^3-3\zeta^2+1\bigr)
+2A_{44}^2(1-\zeta)^3\bigr]
\frac{h_*\Delta_\eta p}{12A_{33}^2 A_{44}}
\nonumber \\
{ } & { } & {}+\frac{A_{13}+A_{44}}{2A_{33}}(1-\zeta)^2
\nabla_\eta\!\cdot{\bf V}^1(\mbox{\boldmath $\eta$})
+C_2(\mbox{\boldmath $\eta$})(1-\zeta),
\label{44ch(4.32)}
\end{eqnarray}
where $C_2(\mbox{\boldmath $\eta$})$ is an arbitrary function.

The substitution of (\ref{44ch(4.32)}) into the boundary condition $(\ref{44ch(4.31)})_1$ yields
$$
C_2=\frac{A_{44}}{2A_{33}^2}h_*\Delta_\eta p
-\frac{A_{44}}{A_{33}}\nabla_\eta\!\cdot{\bf V}^1.
$$

Hence, in view of this relation, formula (\ref{44ch(4.32)}) implies
\begin{equation}
w^2\bigr\vert_{\zeta=0}=-\bigl[(A_{13}+A_{44})^2-4A_{44}^2\bigr]
\frac{h_*\Delta_\eta p}{12A_{33}^2 A_{44}}
+\frac{A_{13}-A_{44}}{2A_{33}}\nabla_\eta\!\cdot{\bf V}^1,
\label{44ch(4.33)}
\end{equation}
where ${\bf V}^1$ is the solution of Eq.~(\ref{44ch(4.29)}).

\section{Local indentation of the coated elastic layer: Leading-order asymptotics for the compressible and incompressible cases}
\label{sec:44.1.4}

Recall that the local indentation of an elastic layer is defined as
$$
w_0({\bf y})\equiv w({\bf y},0),
$$
where $w({\bf y},0)$ is the normal displacement of the layer surface.

In the case of compressible layer, Eqs.~(\ref{44ch(4.24)}) and (\ref{44ch(4.25)}) yield
\begin{equation}
w_0({\bf y})\simeq \frac{h}{A_{33}}p({\bf y}),
\label{44ch(4.s1)}
\end{equation}
so that the deformation response of the coated elastic layer is analogous to that of a Winkler elastic foundation \cite{Johnson1985,PopovHess2013} with the foundation modulus $k=A_{33}/h$. In other words, the deformation of the elastic coating does not contribute substantially to the deformation of a thin compressible layer.

When the layer material approaches the incompressible limit, the right-hand side of (\ref{44ch(4.s1)}) decreases to zero and the first term in the asymptotic expansion (\ref{44ch(4.24)}) disappears. At that, based on the known results \cite{ItskovAksel2002}, it can be shown that the ratios $A_{13}/A_{33}$ and $A_{44}/A_{33}$ tend to 1 and 0, respectively.

Therefore, in the incompressible limit situation formula (\ref{44ch(4.33)}) reduces to
\begin{equation}
w^2\bigr\vert_{\zeta=0}=-\frac{h_*}{12 a_{44}}\Delta_\eta p(\mbox{\boldmath $\eta$})
+\frac{1}{2}\nabla_\eta\!\cdot{\bf V}^1(\mbox{\boldmath $\eta$}),
\label{44ch(4.s2)}
\end{equation}
where $a_{44}=A_{44}$ is the out-of-plane shear modulus of the elastic layer, while
${\bf V}^1(\mbox{\boldmath $\eta$})$ satisfies the equation
\begin{equation}
\frac{1}{h_*}\hat{\mathfrak{L}}^*(\nabla_\eta){\bf V}^1(\mbox{\boldmath $\eta$})
-a_{44}{\bf V}^1(\mbox{\boldmath $\eta$})=\frac{h_*}{2}\nabla_\eta p(\mbox{\boldmath $\eta$}),\quad
\mbox{\boldmath $\eta$}\in \mathbb{R}^2.
\label{44ch(4.s3)}
\end{equation}

Thus, in the case of incompressible bonded elastic layer, formulas  (\ref{44ch(4.2L)})--(\ref{44ch(4.24)}), (\ref{44ch(4.s2)}), and (\ref{44ch(4.s3)}) give
\begin{equation}
w_0({\bf y})\simeq -\frac{h^3}{12 a_{44}}\Delta_y p({\bf y})
+\frac{h}{2}\nabla_y\!\cdot{\bf v}_0({\bf y}),
\label{44ch(4.s4)}
\end{equation}
where the vector ${\bf v}_0({\bf y})$ solves the equation
\begin{equation}
h\hat{\mathfrak{L}}(\nabla_y){\bf v}_0({\bf y})
-a_{44}{\bf v}_0({\bf y})=\frac{h^2}{2}\nabla_y p({\bf y}),\quad
{\bf y}\in \mathbb{R}^2.
\label{44ch(4.s5)}
\end{equation}

Observe that in view of (\ref{44ch(4.28)}), the vector-function ${\bf v}_0({\bf y})$ has the meaning of the tangential displacement of the surface point
$({\bf y},0)$ of the elastic layer.
Recall also that the matrix differential operator $\hat{\mathfrak{L}}(\nabla_y)$ is defined by formulas (\ref{44ch(4.14L)}).

\section{Discussion and conclusion}
\label{sec:44.1.4D}

Observe that in the axisymmetric case (in view of (\ref{44ch(4.13)})), we have
$$
\hat{\bar{\sigma}}_{rr}=
\hat{\bar{A}}_{11} \varepsilon_{rr}
+\hat{\bar{A}}_{12} \varepsilon_{\theta\theta},\quad
\hat{\bar{\sigma}}_{\theta\theta}=
\hat{\bar{A}}_{12} \varepsilon_{rr}
+\hat{\bar{A}}_{11} \varepsilon_{\theta\theta},\quad
\hat{\bar{\sigma}}_{r\theta}=0,
$$
where
$$
\varepsilon_{rr}=\frac{\partial v_r}{\partial r},\quad
\varepsilon_{\theta\theta}=\frac{v_r}{r},
$$
while Eqs.~(\ref{44ch(4.7)}) should be replaced with the following one:
$$
\frac{\hat{h}}{r}\biggl(
\frac{\partial(r\hat{\bar{\sigma}}_{rr})}{\partial r}
-\hat{\bar{\sigma}}_{\theta\theta}\biggr)=-\sigma_{zr}\bigr\vert_{z=0}.
$$

Correspondingly, the boundary condition (\ref{44ch(4.14)}) takes the following form:
\begin{equation}
\sigma_{zr}\bigr\vert_{z=0}=-\hat{h}\hat{\bar{A}}_{11}
\biggl(
\frac{\partial^2 v_r}{\partial r^2}
+\frac{1}{r}\frac{\partial v_r}{\partial r}-\frac{v_r}{r^2}\biggr).
\label{44ch(4.15)}
\end{equation}

It is to note here that the axisymmetric boundary condition (\ref{44ch(4.15)}) was previously derived in a number of papers \cite{AlexandrovMkhitaryan1985,Avilkin1985,RahmanNewaz1997,RahmanNewaz2000}.

Finally, let us consider two opposite limit situations in the incompressible case. First, when the coating is absent and $\hat{\mathfrak{L}}(\nabla_y)\equiv 0$, Eq.~(\ref{44ch(4.s5)}) implies
$$
{\bf v}_0({\bf y})=-\frac{h^2}{2a_{44}}\nabla_y p({\bf y}).
$$
The substitution of this expression into formula (\ref{44ch(4.s4)}) gives
\begin{equation}
w_0({\bf y})\simeq -\frac{h^3}{3 a_{44}}\Delta_y p({\bf y}),
\label{44ch(4.s6)}
\end{equation}
which completely agrees with the asymptotic model developed in \cite{Argatov2012multibody-3,ArgatovMishuris2011b,ArgatovMishuris2011} (see also \cite{Barber1990,Jaffar1989}).

Second, in the case of very stiff (inextensible) coating we will have ${\bf v}_0({\bf y})\equiv{\bf 0}$, and formula (\ref{44ch(4.s4)}) reduces to
\begin{equation}
w_0({\bf y})\simeq -\frac{h^3}{12 a_{44}}\Delta_y p({\bf y}).
\label{44ch(4.s7)}
\end{equation}
In other words, comparing (\ref{44ch(4.s6)}) and (\ref{44ch(4.s7)}),  we conclude that the inextensible membrane coating attached to the surface of a thin incompressible elastic layer reduces the out-of-plane shear compliance of the bonded layer by four times.

The main result of the paper, presented in Section~\ref{sec:44.1.3}, is the leading-order asymptotic formulas for the normal displacement (called the local indentation) of the surface points for compressible and incompressible thin bonded elastic layers reinforced with an elastic membrane.

\section{Acknowledgment}

The authors are grateful for support from the FP7 IRSES Marie Curie grant TAMER No~610547.

\end{document}